\documentclass[reqno, 11pt]{amsart}

\usepackage{amsmath,amssymb,amsthm}
\usepackage[T1]{fontenc}
\usepackage{lmodern}
\usepackage[linktocpage,colorlinks,linkcolor=magenta,citecolor=magenta,urlcolor=magenta,pagebackref=false]{hyperref}
\usepackage{enumitem}
\usepackage[all]{xy}
\usepackage{microtype}
\usepackage{xcolor}
\usepackage{cleveref}  
\usepackage[left=1.25in,right=1.25in,top=1in,bottom=1in]{geometry}

\let\oldphi\phi
\let\phi\varphi
\let\varphi\oldphi
\let\oldepsilon\epsilon
\let\epsilon\varepsilon
\let\varepsilon\oldepsilon

\newcommand{\mbb}[1]{\ensuremath{\mathbb{#1}}}
\newcommand{\Z}{\mbb{Z}}

\newcommand{\mbf}[1]{\ensuremath{\mathbf{#1}}}

\newcommand{\Pb}{\mbf{P}}
\newcommand{\Fb}{\mbf{F}}
\newcommand{\Gb}{\mbf{G}}

\newcommand{\Mb}{\mbf{M}}
\newcommand{\Nb}{\mbf{N}}
\newcommand{\Ab}{\mbf{A}}

\newcommand{\tn}[1]{\textnormal{#1}}
\let\hom\relax\newcommand{\hom}{\tn{Hom}}
\let\ker\relax\newcommand{\ker}{\tn{ker}}
\newcommand{\im}{\tn{im}}

\newcommand{\id}{\tn{id}}
\newcommand{\rank}{\tn{rank}}

\newcommand{\lm}{\tn{lm}}

\newtheorem{thm}{\bf Theorem}[section]

\theoremstyle{definition}
\newtheorem{defn}[thm]{\bf Definition}

\newtheorem{rem}[thm]{\bf Remark}
\newtheorem{ex}[thm]{\bf Example}

\newcommand{\halfskip}{\vspace{0.5em}}

\newcommand{\la}{\langle}
\newcommand{\ra}{\rangle}

\newcommand{\OI}{\tn{OI}}

\numberwithin{equation}{section}
\title{The OIGroebnerBases Package for Macaulay2}

\author{Michael Morrow}
\address{Department of Mathematics, University of Kentucky, 715 Patterson Office Tower, Lexington,
KY 40506 USA}
\email{michael.morrow@uky.edu}

\begin{document}
\begin{abstract}
We introduce the {\it Macaulay2} package {\tt OIGroebnerBases} for working with OI-modules over Noetherian polynomial OI-algebras. The main methods implement OI-analogues of Buchberger's algorithm and Schreyer's theorem to compute Gröbner bases, syzygies and free resolutions of submodules of free OI-modules.
\end{abstract}
\maketitle


\section{Introduction}
\label{section:intro}
Suppose we are given a sequence $(M_n)_{n\in\Z_{\geq0}}$ of related modules $M_n$ over related polynomial rings whose number of variables increases with $n$. One may ask how to simultaneously compute a finite Gröbner basis for each $M_n$. Furthermore, one may ask how to simultaneously compute the module of syzygies of each $M_n$. Using the framework of OI-modules over OI-algebras introduced in \cite{NR19}, these questions were addressed in \cite{MN}, where OI-analogues of Buchberger's algorithm for computing Gröbner bases and Schreyer's theorem for computing syzygies were given. Here, OI denotes the category of totally ordered finite sets and order-preserving increasing maps.

It was further shown in \cite{MN} that the OI-analogue of Schreyer's theorem can be iterated to compute free resolutions of OI-modules out to desired homological degree. Only a few explicit constructions for such resolutions are known; see \cite{FN21,FN22} for examples.

This note introduces the package {\tt OIGroebnerBases}\footnote{Available at \url{https://github.com/morrowmh/OIGroebnerBases}.} for {\it Macaulay2} \cite{GS} to facilitate the computations described above. We review the necessary mathematical background material in \Cref{section:prelim} and summarize the main features of our package in \Cref{section:package}.

\section{Preliminaries}
\label{section:prelim}

We fix notation and recall the needed background on OI-modules. For the rest of this paper, $K$ denotes an arbitrary field.

\begin{defn}
Let OI be the category whose objects are intervals $[n]:=\{1,\ldots,n\}$ for $n\in\Z_{\geq0}$ (we put $[0]=\emptyset$) and whose morphisms are strictly increasing maps $[m]\to[n]$.
\end{defn}
If $\Ab$ is a functor out of OI, we write $\Ab_n$ instead of $\Ab([n])$. We call $\Ab_n$ the \emph{width $n$ component of $\Ab$}. We abuse notation and write $\hom(m,n)$ for the set of all OI-maps $\hom_{\OI}([m],[n])$ from $[m]$ to $[n]$. If $\epsilon\in\hom(m,n)$, we sometimes write $\epsilon_*$ in place of $\Ab(\epsilon)$.
\begin{defn}
\label{def:p}
Let $c>0$ be an integer and define a functor $\Pb=\Pb^c$ from OI to the category of associative, commutative, unital $K$-algebras as follows. For $n\geq0$, define
\[
\Pb_n = K\begin{bmatrix}
x_{1,1}  & \cdots & x_{1,n}\\
\vdots & \ddots & \vdots \\
x_{c,1} & \cdots & x_{c,n}
\end{bmatrix}
\]
and for $\epsilon\in\hom(m,n)$ define $\epsilon_*:\Pb_m\to\Pb_n$ via $x_{i,j}\mapsto x_{i,\epsilon(j)}$.
\end{defn}
Assigning each variable degree $1$, the functor $\Pb$ is a \emph{graded Noetherian polynomial $\OI$-algebra}; see \cite{NR19,MN}.
\begin{defn}[\cite{NR19}]
An OI-\emph{module} $\Mb$ over $\Pb$ is a (covariant) functor from OI to the category of $K$-vector spaces such that
\begin{enumerate}
\item each $\Mb_n$ is an $\Pb_n$-module, and
\item for each $a\in\Pb_m$ and $\epsilon\in\hom(m,n)$ we have a commuting diagram
\[
\label{diagram:oimodule}
\xy\xymatrixrowsep{10mm}\xymatrixcolsep{10mm}
\xymatrix {
\Mb_m\ar[d]_{a\cdot}\ar[r]^{\Mb(\epsilon)} & \Mb_n\ar[d]^{\Pb(\epsilon)(a)\cdot}\\
\Mb_m\ar[r]^{\Mb(\epsilon)} & \Mb_n
}
\endxy
\]
where the vertical maps are multiplication by the indicated elements.
\end{enumerate}
We sometimes refer to $\Mb$ as a $\Pb$-\emph{module}.
\end{defn}
A \emph{homomorphism} of $\Pb$-modules is a natural transformation $\phi:\Mb\to\Nb$ such that each $\phi_n:\Mb_n\to\Nb_n$ is a $\Pb_n$-module homomorphism. We sometimes call $\phi$ a $\Pb$-\emph{linear map}. OI-modules over $\Pb$ and $\Pb$-linear maps form an abelian category with all concepts such as subobject, quotient object, kernel, cokernel, injection, and surjection being defined ``width-wise'' from the corresponding concepts in $K$-vector spaces (see \cite[A.3.3]{W}). Thus, for example, if $\phi:\Mb\to\Nb$ is a $\Pb$-linear map, then the kernel of $\phi$ is a submodule of $\Mb$ defined by $(\ker(\phi))_n=\ker(\phi_n)$. The image of $\phi$ is a submodule of $\Nb$ defined in an analogous fashion.

If $f\in\Mb_n$ for some $n\geq0$ then we call $f$ an \emph{element} of $\Mb$ and write $f\in\Mb$. In this case we say $f$ \emph{has (or is in) width $n$}.
A \emph{subset} of $\Mb$, denoted $S\subseteq\Mb$, is a subset of the disjoint union $\coprod_{n\geq0}\Mb_n$. The submodule of $\Mb$ \emph{generated} by a subset $S\subseteq\Mb$ is the smallest submodule of $\Mb$ containing $S$. This submodule is denoted $\la S\ra_{\Mb}$.

We now discuss freeness.
\begin{defn}[\cite{NR19}]
For any integer $d\geq0$, define an OI-module $\Fb^{\OI,d}$ over $\Pb$ as follows. For $n\in\Z_{\geq0}$ let
\[
\Fb^{\OI,d}_n=\bigoplus_{\pi\in\hom(d,n)}\Pb_n e_\pi\cong (\Pb_n)^{\binom{n}{d}}.  
\]
For $\epsilon\in\hom(m,n)$, define $\Fb^{\OI,d}(\epsilon) \colon \Fb^{\OI,d}_m\to\Fb^{\OI,d}_n$ via $e_\pi\mapsto e_{\epsilon\circ\pi}$. An OI-module $\Fb$ that is isomorphic to a direct sum $\bigoplus_{i=1}^s\Fb^{\OI,d_i}$ for integers $d_1,\ldots,d_s\geq0$ is called a \emph{free} OI-module over $\Pb$ \emph{of rank $s$ generated in widths $d_1,\ldots,d_s$}.
\end{defn}
It is straightforward to see that, given a free OI-module $\Fb=\bigoplus_{i=1}^s\Fb^{\OI,d_i}$, we have
\[
\Fb_n=\bigoplus_{\substack{\pi\in\hom(d_i,n)\\1\leq i\leq s}}\Pb_ne_{\pi,i}
\]
for all $n\geq0$, where the second index on $e_{\pi,i}$ is used to keep track of which direct summand it lives in. We call the $e_{\id_{[d_i]},i}$ the \emph{basis elements} of $\Fb$. The $e_{\id_{[d_i]},i}$ generate $\Fb$ as an OI-module, and to define a $\Pb$-linear map out of $\Fb$ it is enough to specify where the basis elements are mapped. The functor $\Fb$ is an example of a \emph{graded} OI-module \cite{NR19} over $\Pb$ by assigning each basis element degree $0$.

It is convenient to adjust the grading of an OI-module as follows. Given a graded OI-module $\Mb$, define the $d^{th}$ \emph{twist} of $\Mb$ to be the OI-module $\Mb(d)$ that is isomorphic to $\Mb$ as an OI-module, and whose grading is determined by
\[
[\Mb(d)_n]_j=[\Mb_n]_{d+j}.
\]
\begin{ex}
Let $\Pb=\Pb^1$ so that $\Pb_n=K[x_1,\ldots,x_n]$ for $n\geq0$. Then $\Fb^{\OI,1}\oplus\Fb^{\OI,2}$ has its basis elements in degree $0$, while $\Fb^{\OI,1}(-3)\oplus\Fb^{\OI,2}(-4)$ has its basis elements in degrees $3$ and $4$. In width $n$, the rank of both modules as a free $\Pb_n$-module is $\binom{n}{1}+\binom{n}{2}=\binom{n+1}{2}$.
\end{ex}
Let $\Fb=\bigoplus_{i=1}^s\Fb^{\OI,d_i}$ be a free OI-module over $\Pb$ with basis $\{e_{\id_{[d_i]},i}\;:\;i\in[s]\}$. A \emph{monomial} in $\Fb$ is an element of the form $ae_{\pi,i}$ where $a$ is a monomial in $\Pb$. There is a suitable notion of a \emph{monomial order} on the monomials of $\Fb$ (see \cite[Definition 3.1 and Example 3.2]{MN}) with which we can define the \emph{lead monomial} $\lm(f)$ of any element $f\in\Fb$. Moreover, we define $\lm(E)=\{\lm(f)\;:\;f\in E\}$ for any subset $E\subseteq \Fb$.

Our primary object of study is defined as follows.
\begin{defn}[\cite{NR19,MN}]
Fix a monomial order $<$ on $\Fb$ and let $\Mb$ be a submodule of $\Fb$. A subset $G\subseteq\Mb$ is called a \emph{Gröbner basis} of $\Mb$ (with respect to $<$) if
\[
\la\lm(\Mb)\ra_{\Fb}=\la\lm(G)\ra_{\Fb}.
\]
\end{defn}
In \cite{NR19}, it was established that any submodule of a finitely generated free OI-module over a Noetherian polynomial OI-algebra has a finite Gröbner basis. It was shown in \cite{MN} how to compute such bases in finite time. Our package implements this construction with the {\tt oiGB} method; see \Cref{subs:oigb}.

We also consider syzygies. Given a finitely generated submodule $\Mb$ of a free OI-module $\Fb$, there is a canonical surjective $\Pb$-linear map $\phi:\Gb\to\Mb$ sending the basis elements of a free OI-module $\Gb$ to the generators of $\Mb$ (see \cite[Proposition 3.19]{NR19}). The {\tt oiSyz} method in our package implements a construction given in \cite{MN} for computing the kernel of $\phi$; see \Cref{subs:oisyz}.

Finally, it was shown in \cite{MN} how to iterate the syzygy construction to compute free resolutions $\Fb^\bullet\to\Mb\to0$ out to desired homological degree. If $\Mb$ is graded, then $\Fb^\bullet$ can be pruned in order to form a \emph{graded minimal free resolution} of $\Mb$ (see \cite{FN21} and \cite[Theorem 5.4]{MN}). This is implemented in our package with the {\tt oiRes} method; see \Cref{subs:oires}.

\section{The Package}
\label{section:package}
The main methods of our package are {\tt oiGB} for computing Gr\"obner bases, {\tt oiSyz} for computing syzygies and {\tt oiRes} for computing resolutions. This section illustrates how to use these methods. For more information about other methods, as well as optional arguments such as grading shifts, we refer the reader to the package documentation.
\subsection{Gr\"obner bases}
\label{subs:oigb} Let $\Fb$ be a finitely generated free $\Pb$-module and let $\Mb$ be a submodule generated by $B=\{b_1,\ldots,b_s\}$. Fix a monomial order $<$ on $\Fb$. Then by \cite[Algorithm 3.17]{MN}, a finite Gröbner basis $G$ (with respect to $<$) for $\Mb$ containing $B$ can be computed in finite time. Using our package, one computes such Gröbner bases with the {\tt oiGB} method.

\begin{ex}
Let $\Fb=\Fb^{\OI,1}\oplus\Fb^{\OI,1}\oplus\Fb^{\OI,2}$ have basis $\{e_{\id_{[1]},1},e_{\id_{[1]},2},e_{\id_{[2]},3}\}$. Let $\Pb=\Pb^2$ so that $\Pb$ has two rows of variables, and let
\[
B=\{x_{1,1}e_{\id_{[1]},1}+x_{2,1}e_{\id_{[1]},2},\;x_{1,2}x_{1,1}e_{\pi,2}+x_{2,2}x_{2,1}e_{\id_{[2]},3}\}
\]
where $\pi:[1]\to[2]$ is given by $1\mapsto 2$. Thus, the first element of $B$ has width $1$ and the second element has width $2$. Fix the \emph{lex order} on $\Fb$ as described in \cite[Example 3.2]{MN}. We compute a finite Gröbner basis for $\la B\ra_{\Fb}$ in {\it Macaulay2} as follows. First, we define our polynomial OI-algebra $\Pb$ with {\tt makePolynomialOIAlgebra}. The user must specify the number of variable rows, the variable symbol, and the ground field $K$:
\halfskip
\begin{verbatim}
i1 : needsPackage "OIGroebnerBases";
i2 : P = makePolynomialOIAlgebra(2, x, QQ);
\end{verbatim}
\halfskip
Now we define our free OI-module $\Fb$ with {\tt makeFreeOIModule}. The user specifies the basis symbol, a list of basis element widths, and the underlying polynomial OI-algebra:
\halfskip
\begin{verbatim}
i3 : F = makeFreeOIModule(e, {1,1,2}, P);
\end{verbatim}
\halfskip
Since we want to define our elements of $B$, we need to call the {\tt installGeneratorsInWidth} method in order to work with our basis symbol {\tt e}. This method takes a free OI-module and a width as input:
\halfskip
\begin{verbatim}
i4 : installGeneratorsInWidth(F, 1);
i5 : installGeneratorsInWidth(F, 2);
\end{verbatim}
\halfskip
We're ready to define the elements of $B$:
\halfskip
\begin{verbatim}
i6 : use F_1; b1 = x_(1,1)*e_(1,{1},1)+x_(2,1)*e_(1,{1},2);
i8 : use F_2; b2 = x_(1,2)*x_(1,1)*e_(2,{2},2)+x_(2,2)*x_(2,1)*e_(2,{1,2},3);
\end{verbatim}
\halfskip
Here, for example, {\tt e\_(2,\{2\},2)} is the element $e_{\pi,2}$ as defined above. In general, an element $e_{\sigma,i}\in\Fb$ translates to an object in our package as follows. Suppose $\sigma\in\hom(m,n)$, and write $\im(\sigma)=\{a_1,\ldots,a_m\}$ where each $a_j\in[n]$. Then $e_{\sigma,i}$ becomes {\tt e\_(n,\{a\textsubscript{1},\ldots,a\textsubscript{m}\},i)} in our package. Now let's compute a Gröbner basis with the method {\tt oiGB}, which takes a list of elements as input:
\halfskip
\begin{verbatim}
i10 : oiGB {b1, b2}
o10 = {x   e        + x   e       , x   x   e        + x   x   e          ,
        1,1 1,{1},1    2,1 1,{1},2   1,2 1,1 2,{2},2    2,2 2,1 2,{1, 2},3
      ---------------------------------------------------------------------
      x   x   x   e           - x   x   x   e          }
       2,3 2,2 1,1 3,{2, 3},3    2,3 2,1 1,2 3,{1, 3},3
\end{verbatim}
\halfskip
This tells us that a Gröbner basis for $\Mb=\la B\ra_{\Fb}$ with respect to the lex order is given by the following elements:
\begin{align*}
b_1=x_{1,1}e_{\id_{[1]},1}+x_{2,1}e_{\id_{[1]},2}&\in\Fb_1\\
b_2=x_{1,2}x_{1,1}e_{\pi,2}+x_{2,2}x_{2,1}e_{\id_{[2]},3}&\in\Fb_2\\
b_3=x_{2,3}x_{2,2}x_{1,1}e_{\sigma_1,3}-x_{2,3}x_{2,1}x_{1,2}e_{\sigma_2,3}&\in\Fb_3
\end{align*}
where $\sigma_1:[2]\to[3]$ is given by $1\mapsto 2$ and $2\mapsto 3$ and $\sigma_2:[2]\to[3]$ is given by $1\mapsto 1$ and $2\mapsto 3$. This agrees with \cite[Example 3.20]{MN}. It follows that, given any $n\geq3$, a finite Gröbner basis for $\Mb_n$ is given by the images of $b_1$, $b_2$ and $b_3$ under any morphism $[1]\to[n]$, $[2]\to[n]$ and $[3]\to[n]$ respectively.
\end{ex}

\begin{rem}
Passing the optional argument {\tt Verbose => true} into the methods {\tt oiGB}, {\tt oiSyz} and {\tt oiRes} will print useful debug information, and also provides a way to track the progress of the computation.
\end{rem}

\subsection{Syzygies}
\label{subs:oisyz} As in the previous section, let $\Fb$ be a finitely generated free $\Pb$-module. Let $B=\{b_1,\ldots,b_s\}\subset\Fb$ and let $w_i$ denote the width of $b_i$. Define the free OI-module $\Gb=\bigoplus_{i=1}^s\Fb^{\OI,w_i}$ with basis $\{d_{\id_{[w_i]},i}\}$. Let $\phi:\Gb\to\la B\ra_{\Fb}$ be the canonical surjective map defined by $d_{\id_{[w_i]},i}\mapsto b_i$. The {\it syzygy module} of $\la B\ra_{\Fb}$ is defined to be the kernel of $\phi$. Suppose $B$ is a Gröbner basis for $\la B\ra_{\Fb}$ with respect to some monomial order $<$ on $\Fb$. Using the construction described in \cite[Theorem 4.6]{MN}, the method {\tt oiSyz} computes a finite Gröbner basis for $\ker(\phi)$ with respect to a suitable monomial order on $\Gb$ induced by $<$.

\begin{ex}
\label{ex:syz}
Let $\Pb=\Pb^2$ and let $\Fb=\Fb^{\OI,1}\oplus\Fb^{\OI,1}$ have basis $\{e_{\id_{[1]},1},e_{\id_{[1]},2}\}$. Define
\[
f=x_{1,2}x_{1,1}e_{\pi,1}+x_{2,2}x_{2,1}e_{\rho,2}\in\Fb_2
\]
where $\pi:[1]\to[2]$ is given by $1\mapsto2$ and $\rho:[1]\to[2]$ is given by $1\mapsto1$. We will compute a Gröbner basis $G$ for $\la f\ra_{\Fb}$, and then compute the syzygy module of $G$. Starting a new {\it Macaulay2} session, we run the following:
\halfskip
\begin{verbatim}
i1 : needsPackage "OIGroebnerBases";
i2 : P = makePolynomialOIAlgebra(2, x, QQ);
i3 : F = makeFreeOIModule(e, {1,1}, P);
i4 : installGeneratorsInWidth(F, 2);
i5 : use F_2; f = x_(1,2)*x_(1,1)*e_(2,{2},1)+x_(2,2)*x_(2,1)*e_(2,{1},2);
i7 : G = oiGB {f}
o7 = {x   x   e        + x   x   e       ,
       1,2 1,1 2,{2},1    2,2 2,1 2,{1},2
     --------------------------------------------
     x   x   x   e        - x   x   x   e       }
      2,3 2,2 1,1 3,{2},2    2,3 2,1 1,2 3,{1},2
\end{verbatim}
\halfskip
Hence, $\la f\ra_{\Fb}$ has a Gröbner basis (with respect to the lex order)
\[
G =\{x_{1,2}x_{1,1}e_{\pi,1}+x_{2,2}x_{2,1}e_{\rho,2},\;x_{2,3}x_{2,2}x_{1,1}e_{\sigma_1,2}-x_{2,3}x_{2,1}x_{1,2}e_{\sigma_2,2}\}
\]
where $\sigma_1:[1]\to[3]$ is given by $1\mapsto2$ and $\sigma_2:[1]\to[3]$ is given by $1\mapsto1$. Define the free OI-module $\Gb=\Fb^{\OI,2}(-2)\oplus\Fb^{\OI,3}(-3)$ with basis $\{d_{\id_{[2]},1},d_{\id_{[3]},2}\}$. The package assigns these degree shifts automatically. Putting $g=x_{2,3}x_{2,2}x_{1,1}e_{\sigma_1,2}-x_{2,3}x_{2,1}x_{1,2}e_{\sigma_2,2}\in\Gb_3$ so that $G=\{f,g\}$, we define the map $\phi:\Gb\to\la G\ra_{\Fb}$ via $d_{\id_{[2]},1}\mapsto f$ and $d_{\id_{[3]},2}\mapsto g$. We can now compute a Gröbner basis $D$ for $\ker(\phi)$ (with respect to the Schreyer order on $\Gb$ induced by the monomial order on $\Fb$; see \cite[Definition 4.2]{MN}) using the method {\tt oiSyz}. The user inputs the Gröbner basis $G$ and the basis symbol $d$:
\halfskip
\begin{verbatim}
i8 : D = oiSyz(G, d)
o8 = {x   d           - x   d           + 1d             ,
       1,2 3,{1, 3},1    1,1 3,{2, 3},1     3,{1, 2, 3},2
     -----------------------------------------------------------
     x   d              - x   d             , x   d
      2,4 4,{1, 2, 3},2    2,3 4,{1, 2, 4},2   1,2 4,{1, 3, 4},2
     -----------------------------------------------------------
     - x   d              - x   d             }
        1,1 4,{2, 3, 4},2    1,3 4,{1, 2, 4},2
\end{verbatim}
\halfskip
This says that $\ker(\phi)$ has a Gröbner basis $D$ given by the following elements of $\Gb$:
\begin{align*}
x_{1,2}d_{\pi_1,1}-x_{1,1}d_{\pi_2,1}-d_{\id_{[3]},2}&\in\Gb_3\\
x_{2,4}d_{\pi_3,2}-x_{2,3}d_{\pi_4,2}&\in\Gb_4\\
x_{1,2}d_{\pi_5,2}-x_{1,1}d_{\pi_6,2}-x_{1,3}d_{\pi_4,2}&\in\Gb_4
\end{align*}
where the $\pi_i$ for $1\leq i\leq 6$ are given as follows:
\begin{align*}
&\pi_1:[2]\to[3]\quad\text{via}\quad1\mapsto1,2\mapsto3\\
&\pi_2:[2]\to[3]\quad\text{via}\quad1\mapsto2,2\mapsto3\\
&\pi_3:[3]\to[4]\quad\text{via}\quad1\mapsto1,2\mapsto2,3\mapsto3\\
&\pi_4:[3]\to[4]\quad\text{via}\quad1\mapsto1,2\mapsto2,3\mapsto4\\
&\pi_5:[3]\to[4]\quad\text{via}\quad1\mapsto1,2\mapsto3,3\mapsto4\\
&\pi_6:[3]\to[4]\quad\text{via}\quad1\mapsto2,2\mapsto3,3\mapsto4.
\end{align*}
\end{ex}

\subsection{Resolutions}
\label{subs:oires} The syzygy construction described in \cite[Theorem 4.6]{MN} can be iterated to build resolutions of submodules of free OI-modules. Let $\Fb$ be a free $\Pb$-module of finite rank, and let $\Mb\subseteq\Fb$ be a submodule generated by a finite set $B$. Then a free resolution of $\Mb$ can be computed out to desired homological degree using \cite[Procedure 5.1]{MN}. Moreover, if $\Mb$ is homogeneous, then a graded minimal free resolution of $\Mb$ can be computed out to arbitrary homological degree.
\begin{ex}
Let $\Pb=\Pb^2$ and let $\Fb=\Fb^{\OI,1}\oplus\Fb^{\OI,1}$ have basis $\{e_{\id_{[1]},1},e_{\id_{[1]},2}\}$, so $\Fb$ has rank 2. Define
\[
f=x_{1,2}x_{1,1}e_{\pi,1}+x_{2,2}x_{2,1}e_{\rho,2}\in\Fb_3
\]
where $\pi:[1]\to[3]$ is given by $1\mapsto2$ and $\rho:[1]\to[3]$ is given by $1\mapsto1$. Note the similarity to \Cref{ex:syz}. Since $f$ is homogeneous, $\la f\ra_{\Fb}$ is a graded submodule, and we will compute the beginning of a graded minimal free resolution using {\tt oiRes}. The user specifies a list of elements (who generate the module to be resolved) and a homological degree. In a new {\it Macaulay2} session, we run the following:
\halfskip
\begin{verbatim}
i1 : needsPackage "OIGroebnerBases";
i2 : P = makePolynomialOIAlgebra(2, x, QQ);
i3 : F = makeFreeOIModule(e, {1, 1}, P);
i4 : installBasisElements(F, 3);
i5 : use F_3; f = x_(1,2)*x_(1,1)*e_(3,{2},1)+x_(2,2)*x_(2,1)*e_(3,{1},2);
i7 : ranks oiRes({f}, 5)
o7 = 0: rank 1
     1: rank 2
     2: rank 4
     3: rank 7
     4: rank 11
     5: rank 22
\end{verbatim}
\halfskip
Note: if one computes out to homological degree $n$, then only the first $n-1$ ranks are guaranteed to be minimal. Thus, we have the beginning of a minimal free resolution for $\Mb=\la f\ra_{\Fb}$:
\[
\cdots \to \Fb^4\to\Fb^3\to\Fb^2\to\Fb^1\to\Fb^0\to\Mb\to 0
\]
where
\begin{align*}
\rank(\Fb^0)&=1\\
\rank(\Fb^1)&=2\\
\rank(\Fb^2)&=4\\
\rank(\Fb^3)&=7\\
\rank(\Fb^4)&=11.
\end{align*}
One can obtain more information about resolutions such as grading shifts, generators of the free modules, and differentials by using the {\tt describe} method. We refer the reader to the package documentation. Such information can be used to restrict a resolution of $\Mb$ to any width $w$ to obtain a graded (but not necessarily minimal) free resolution of the $\Pb_w$-module $\Mb_w$, as in \cite[Section 3]{FN21}.
\end{ex}


\end{document}